\documentclass{amsart}

\newtheorem{THM}{Theorem}[section]
\newtheorem{LMA}[THM]{Lemma}
\newtheorem{PROP}[THM]{Proposition}
\newtheorem{CORO}[THM]{Corollary}
\newtheorem{CONJ}[THM]{Conjecture}

\numberwithin{equation}{section}

\usepackage{amsmath}
\usepackage{amssymb}

\newcommand{\showon}{\begin{eqnarray}}
\newcommand{\showoff}{\end{eqnarray}}

\newcommand{\NN}{\mathbb{N}}
\newcommand{\RR}{\mathbb{R}}
\newcommand{\ZZ}{\mathbb{Z}}
\newcommand{\CC}{\mathbb{C}}
\newcommand{\R}[1]{R_{#1}}
\newcommand{\HH}[1]{H_{#1}}

\newcommand{\J}[1]{J_{#1}}
\newcommand{\Jo}[1]{J_{0}^{#1}}
\newcommand{\Ji}[1]{J_{1}^{#1}}
\newcommand{\Ko}[1]{K_{0}^{#1}}
\newcommand{\Ki}[1]{K_{1}^{#1}}
\newcommand{\PP}{\mathbf{P}}
\newcommand{\Q}{\mathbf{Q}}
\newcommand{\I}{\mathrm{I}}
\newcommand{\II}{\mathbb{I}}
\newcommand{\al}{\alpha}
\newcommand{\lam}{\lambda}
\newcommand{\drop}{\smallsetminus}
\newcommand{\eps}{\varepsilon}
\newcommand{\bul}{\bullet}
\newcommand{\lraw}{\longrightarrow}
\newcommand{\llaw}{\longleftarrow}
\newcommand{\uaw}{\uparrow}
\newcommand{\SQ}{\boxplus}
\newcommand{\tri}{\triangledown}

\begin{document}

\title[Reliability polynomials]{Zeros of reliability polynomials and
$f$-vectors of matroids}
\author{David G. Wagner}
\address{Department of Combinatorics and Optimization\\
University of Waterloo\\
Waterloo, Ontario, Canada\ \ N2L 3G1}
\email{{\tt dgwagner@math.uwaterloo.ca}}
\thanks{Research supported by the Natural
Sciences and Engineering Research Council of Canada under
operating grant OGP0105392.}
\keywords{reliability polynomials, series-parallel networks,
Schur stability, Hurwitz stability, Hermite-Biehler theorem,
$f$-vectors and $h$-vectors of matroids.}
\subjclass{05C99, 26C10, 95C15, 06A08.}

\begin{abstract}
For a finite multigraph $G$, the {\em reliability function} of $G$
is the probability $\R{G}(q)$ that if each edge of $G$ is deleted
independently with probability $q$ then the remaining edges of
$G$ induce a connected spanning subgraph of $G$;  this is a polynomial
function of $q$.  In 1992,  Brown and Colbourn conjectured that for any
connected multigraph $G$, if $q\in\CC$ is such that $\R{G}(q)=0$ then
$|q|\leq 1$.   We verify that this conjectured property of $\R{G}(q)$
holds if $G$ is a series-parallel network.  The proof is by an
application of the Hermite-Biehler Theorem and development of a
theory of higher-order interlacing for polynomials with only real
nonpositive zeros.  We conclude by establishing some new
inequalities which are satisfied by the $f$-vector of any
matroid without coloops, and by discussing some stronger
inequalities which would follow (in the cographic case) from
the Brown-Colbourn Conjecture, and are hence true for cographic
matroids of series-parallel networks.
\end{abstract}
\maketitle

\section{Introduction}

Given a  finite multigraph $G=(V,E)$ and $0\leq q\leq 1$, let
$\mathcal{G}(q)$ denote a random spanning subgraph of $G$ obtained by
deleting each edge of $G$ independently with probability $q$.
The {\em reliability function} of $G$ is the probability $\R{G}(q)$
that $\mathcal{G}(q)$ is connected, considered as a function of $q$.
(Trivially, if $G$ is not connected then $\R{G}(q)\equiv 0$
identically.)  In fact (as we see in Section 1) this is a polynomial
function of $q$.  In 1992, Brown and Colbourn \cite{BC} made the
following fascinating conjecture.  A polynomial $S(q)$ is {\em
Schur quasi-stable} if every $q\in\CC$ for which $S(q)=0$
is such that $|q|\leq 1$; for the relevance of this concept to
solutions of linear finite difference equations, see Theorem 3.2 of
Barnett \cite{Ba}.
\begin{CONJ}[Brown-Colbourn] For any connected multigraph $G$, the
reliability polynomial $R_{G}(q)$ is Schur quasi-stable.\end{CONJ}
In support of this conjecture, Brown and Colbourn verify that
this property holds for all simple graphs on up to six vertices,
and show that for every multigraph $G$ there is a multigraph $G'$
which is obtained from $G$ by repeatedly subdividing edges, and for
which $R_{G'}(q)$ is Schur quasi-stable.  The proofs in \cite{BC}
are based on the Enestr\"om-Kakeya Theorem, which gives a
sufficient condition for a polynomial with real coefficients to
be Schur quasi-stable.  As Brown and Colbourn remark, however,
there are multigraphs for which the Enestr\"om-Kakeya Theorem fails
to show that $\R{G}(q)$ is Schur quasi-stable.  Some other explanation
must be sought if the Brown-Colbourn Conjecture is to be proven.

In this paper we give some indications that the Hermite-Biehler
Theorem can provide such an explanation.
This theorem is a necessary and sufficient condition for a
polynomial $P(u)$ with real coefficients to be such that all
its zeros have nonpositive real part;  by a fractional linear
transformation we can map the unit disc to left half-plane
and apply the Hermite-Biehler Theorem to reliability polynomials.
Informally, the condition is that if $P(u)$ is expanded into its
even and odd parts, $P(u)=P_{0}(u^{2})+uP_{1}(u^{2})$, then the
polynomials $P_{1}(x)$ and $P_{0}(x)$ have all their zeros on
the nonpositive part of the real axis, and these zeros ``interlace''
(in a sense we make precise in Section 2).  This allows the
well-developed theory of polynomials with only real zeros
to be applied to Conjecture 0.1, but even this is not sufficient.
More precisely, this theory must be developed further
in order to obtain significant results on reliability polynomials.

The key extension of technique in this paper is the introduction of
a useful concept of higher order interlacing for polynomials
with only real nonpositive zeros; this involves the definition of
``interpolatory hypercubes of polynomials'' of any dimension.
However, because of the complexity of the relations derived from
two-vertex-cut reduction for reliability polynomials,
we have applied this theory here for interpolatory cubes only up
to dimension four.  The limited scope of Theorem 0.2 in
relation to the generality of Conjecture 0.1 reflects only this
artificial restriction, and should not be interpreted as stemming
from some intrinsic limitation of the method.

The class $\mathfrak{SP}$ of {\em series-parallel networks} is defined
recursively as follows.  Every multigraph $G$ in $\mathfrak{SP}$ has a
distinguished unordered pair $\{s,t\}$ of distinct vertices, called the
{\em terminals} of $G$.  If $G$ consists of just one edge connecting its
terminals, then $G$ is in $\mathfrak{SP}$.  Let $G$ and $G'$ be in $\mathfrak{SP}$,
with terminals $\{s,t\}$ and $\{s',t'\}$, respectively.  If $G$ and $G'$
have no edges in common, and only the vertex $t=s'$ in common, then
$G\cup G'$ is in $\mathfrak{SP}$, and is called a {\em series connection}
of $G$ and $G'$; its terminals are $\{s,t'\}$.  If $G$ and $G'$ have
no edges in common, and only the vertices $s=s'$ and $t=t'$ in common,
then $G\cup G'$ is in $\mathfrak{SP}$, and is called a
{\em parallel connection} of $G$ and $G'$; its terminals are
$\{s,t\}$.  Let $\mathfrak{SP}'$ denote the class of connected multigraphs
every two-connected component of which is in the class $\mathfrak{SP}$
(for some choice of terminals).

\begin{THM}  If the multigraph $G$ is in the class $\mathfrak{SP}'$
then $R_{G}(q)$ is Schur quasi-stable.\end{THM}

Since at least the early 1970s there has been some interest in
obtaining inequalities valid for the $f$-vectors and/or the
$h$-vectors of simplicial complexes belonging to various classes;
in part, this developed from similar investigations in the 19th century
into the combinatorial geometry of convex polytopes.
Sections II.2, II.3, III.1, and III.3 of Stanley \cite{St2} provide
an excellent overview of these results.  Also, Section 5 of
Bj\"orner \cite{Bj} considers in detail the case of matroids,
and Ball and Provan \cite{BP} and Colbourn \cite{Co} discuss the
application of these ideas to estimation of network reliability.
In contrast with the cases of Cohen-Macaulay complexes or simplicial
polytopes, the case of matroids is still only rather poorly understood;
some recent results in this direction are Brown and Colbourn \cite{BC2}
and Chari \cite{Ch1,Ch2}.  In fact, Conjecture 0.1 implies numerous
strong inequalities for the $f$-vector of a cographic matroid,
as we shall see in Section 6;  thus, by Theorem 0.2, these
inequalities hold for cographic matroids of multigraphs in the class
$\mathfrak{SP}'$.  Moreover, it is an elementary
consequence of Chari's recent work \cite{Ch1,Ch2} that the weakest
of these inequalities hold more generally for matroids.

\begin{THM}  Let $M$ be (the set of independent sets of) a matroid
of rank $d$, and for $0\leq i\leq d$ let $f_{i}$ denote the
number of $i$-element sets in $M$.  If $M$ has no coloops then
for all $0\leq k\leq d$,
$$0\leq \sum_{i=k}^d\binom{i}{k}(-2)^{d-i}f_i.$$
\end{THM}
These inequalities are violated by some simplicial polytopes
and some broken-circuit complexes, and are satisfied with
equality for all $0\leq k<d$ if $M$ is a direct sum of $2$-circuits.
It thus appears that a better understanding of the phenomena
underlying Conjecture 0.1 could lead not only to improved methods for
estimating network reliability, but perhaps toward a set of strong
necessary conditions on the $f$-vectors of matroids in general.

In Section 1 we review the bare essentials of the combinatorics
of reliability polynomials, the deletion/contraction algorithm
and two-vertex-cut reduction, and we translate Conjecture 0.1 into
a form to which the Hermite-Biehler Theorem applies.  (No familiarity
with matroid theory is assumed until Section 6.)  In Section 2 we review the
Hermite-Biehler Theorem and state the lemmas on polynomials with
only real zeros which are useful.  In Section 3 we sketch how just this
amount of theory can be used to verify Conjecture 0.1 for all multigraphs
such that the underlying simple graph of every two-connected component is
an edge or a cycle.  Section 4 contains the new theoretical
development of the paper;  in Section 5 we apply this technique
to prove Theorem 0.2.  We conclude in Section 6 with a discussion
of reliability polynomials in the more general context of
Cohen-Macaulay complexes, for which we assume familiarity with the
standard concepts.  Readers interested specifically in Theorem 0.3
can skip directly to Section 6.

I gratefully acknowledge Jason Brown and Charlie Colbourn for telling
me of their conjecture, J\"urgen Garloff for telling me of the
Hermite-Biehler Theorem, and Bruce Richmond for several stimulating
conversations.

\section{Reliability Polynomials}

For a more thorough introduction to reliability polynomials, see
Colbourn \cite{Co}.  By a {\em multigraph} we mean a  finite
graph which may possess both loops and multiple edges.
It is clear that for multigraphs $G$ and $G'$,
\showon \mbox{if $G\simeq G'$ then}\  \R{G}(q)=\R{G'}(q). \showoff
If $G$ and $N$ are multigraphs with exactly one
vertex in common then
\showon \R{G\cup N}(q)=\R{G}(q)\R{N}(q), \showoff
as follows directly from the definition.
For a  multigraph $G=(V,E)$ and any $e\in E$ let
$G\drop e$ denote $G$ with $e$ deleted and let $G/e$ denote $G$ with
$e$ contracted; then
\showon \R{G}(q)=q\R{G\drop e}(q)+(1-q)\R{G/e}(q), \showoff
since the conditional probability that $\mathcal{G}(q)$ is connected
given that $e$ is deleted is $\R{G\drop e}(q)$, and the conditional
probability that $\mathcal{G}(q)$ is connected given that $e$ is not
deleted is $\R{G/ e}(q)$.
If $G'$ is obtained from $G$ by removing all loops, then
\showon \R{G}(q)=\R{G'}(q), \showoff
since if $e$ is a loop of $G$ then $G\drop e\simeq G/e$, and we can
apply (1.1) and (1.3) and induction on the number of loops of $G$.
Henceforth by a  {\em network}  we shall mean a  finite connected graph
which has no loops but may have multiple edges.

For a network $G$ we denote by $G^{\natural}$ the {\em underlying simple
graph of $G$},  which  has the same vertices as $G$ and one edge $u\sim v$
for every pair of vertices $\{u,v\}$ which are adjacent in $G$.
A {\em spindle} in a network $G$ is a (nonempty) maximal set of edges
in $G$ all of  which are incident with the same pair of (distinct) vertices  of
$G$; if the spindle has $c$ edges then we say it is a {\em $c$-spindle}.
There is an obvious natural bijection between the spindles of $G$ and
the edges of $G^{\natural}$.  If $\sigma$ is a spindle in $G$ then let
$G\drop\sigma$ be obtained from $G$ by deleting all the edges of $\sigma$,
and let $G/\sigma$ be obtained from $G$ by contracting all the edges of
$\sigma$\ (notice that $G/\sigma$ has no loops).  If $\sigma$ is a $c$-spindle
in the network $G$ then
\showon \R{G}(q)=q^{c}\R{G\drop\sigma}(q)+(1-q^{c})\R{G/\sigma}(q), \showoff
as follows from (1.3) and (1.4) by induction on $c$.

As examples, let $kT_{n}$ denote a network for which the
underlying simple graph is a tree with $n$ vertices and in which
each spindle has $k$ edges, and let $kC_{n}$ denote a network for
which the underlying simple graph is a cycle with $n$ vertices and in
which each spindle has $k$ edges.  Since $R_{kT_{2}}(q)=1-q^{k}$, it
follows by (1.2) and induction on $n$ that
$R_{kT_{n}}(q)=(1-q^{k})^{n-1}$ for all $k\geq 1$ and $n\geq 2$.
From this, (1.5), and induction on $n$ it follows that
$\R{kC_{n}}(q)=(1-q^{k})^{n-1}(1+(n-1)q^{k})$ for all $k\geq 1$
and $n\geq 3$.  As noted by Brown and Colbourn \cite{BC}
(Proposition 5.1) these examples suffice to show that the closure
of the set of all zeros of reliability polyomials contains the
whole unit disc $\{q\in\CC:\ |q|\leq 1\}$.

Generalizing (1.3) and (1.5), let $G$ and $N$ be two networks which
intersect in exactly two vertices $v$ and $w$, and let $G^{\bul}$
and $N^{\bul}$ be obtained by identifying $v$ and $w$ in $G$ and in
$N$, respectively, and removing any loops thus produced.  Then
\showon \R{G\cup N}(q)=\R{G}(q)\R{N^{\bul}}(q)+
\R{G^{\bul}}(q)\R{N}(q)-\R{G}(q)\R{N}(q),\showoff
since the conditional probability that $\mathcal{G}(q)\cup\mathcal{N}(q)$
is connected given that $\mathcal{G}(q)$ is connected is
$\R{N^{\bul}}(q)$, the second term has a similar interpretation,
the third term corrects double counting of the case that both
$\mathcal{G}(q)$ and $\mathcal{N}(q)$ are connected, and if neither
$\mathcal{G}(q)$ nor $\mathcal{N}(q)$ is connected then
$\mathcal{G}(q)\cup\mathcal{N}(q)$ is not connected.

For a network $G$ we denote by $m_{G}$ the number of edges
of $G$ and by $n_{G}$ the number of vertices of $G$, and we let
$d_{G}:=m_{G}-n_{G}+1$;\ we omit the subscript when no confusion
can arise.  From (1.5) it follows by induction on $m$ that
$\R{G}(q)$ is a polynomial in $\ZZ[q]$ of degree $m$,
and the multiplicity of $q=1$ as a zero of $\R{G}(q)$ is at least
$n-1$.  In view of this, for each network $G$ we may define the
polynomial
\showon \HH{G}(q):=(1-q)^{1-n}\R{G}(q)\showoff
in $\ZZ[q]$.  The Brown-Colbourn Conjecture is
equivalent to the claim that for any network $G$,
\showon \mbox{if $q\in\CC$ is such that $\HH{G}(q)=0$ then
$|q|\leq 1$.} \showoff
It follows from (1.2) and (1.7) that if $G$ and $N$ are networks
with exactly one vertex in common then
\showon \HH{G\cup N}(q)=\HH{G}(q)\HH{N}(q).\showoff
If $\sigma$ is a $c$-spindle in the network $G$ then
\showon \HH{G}(q)=q^{c}\HH{G\drop\sigma}(q)+\left(\frac{1-q^{c}}{1-q}\right)
\HH{G/\sigma}(q), \showoff  as follows from (1.5) and (1.7).
Similarly, with notation as in (1.6) we see that
\showon \quad\quad\HH{G\cup N}(q)=\HH{G}(q)\HH{N^{\bul}}(q)+
\HH{G^{\bul}}(q)\HH{N}(q)-(1-q)\HH{G}(q)\HH{N}(q).\showoff
It follows from (1.10) by induction on $m$ that
$\HH{G}(q)$ is a polynomial of degree $d$ with nonnegative integer
coefficients, and the constant term of $\HH{G}(q)$ is $1$.  (In fact,
the coefficients of $\HH{G}(q)$ form the {\em $h$-vector} of the
cographic matroid of $G$ and have been studied extensively in
the context of Cohen-Macaulay simplicial complexes.  We shall
return to this point in Section 6.)

We make a  change of variables $u:=(-1-q)/(1-q)$ and
conversely $q=(-1-u)/(1-u)$, and define
\showon \J{G}(u):=(u-1)^{d}\HH{G}\left(\frac{-1-u}{1-u}\right),\showoff
so that
\showon \HH{G}(q)=\left(\frac{q-1}{2}\right)^{d}
\J{G}\left(\frac{-1-q}{1-q}\right). \showoff
From (1.12) and (1.13) it follows that (1.8) is equivalent to
\showon \mbox{if $u\in\CC$ is such that $\J{G}(u)=0$ then
$\Re(u)\leq 0$.} \showoff
From (1.9) and (1.12) it follows that if $G$ and $N$ are networks with
exactly one vertex in common then
\showon \J{G\cup N}(u)=\J{G}(u)\J{N}(u).\showoff
It follows from (1.10) and (1.12) that if $\sigma$ is a $c$-spindle in
the network $G$ then
\showon \J{G}(u) = (u+1)^{c}\J{G\drop\sigma}(u)+
\left[\frac{(u+1)^{c}-(u-1)^{c}}{2}\right]\J{G/\sigma}(u).\showoff
Similarly, with notation as in (1.6) and (1.11) we see that
\showon \J{G\cup N}(u)=\J{G}(u)\left[\J{N}(u)+\J{N^{\bul}}(u)\right]+
\left[\J{G}(u)+\J{G^{\bul}}(u)\right]\J{N}(u).\showoff
From (1.16) and induction on $m$ it follows that $\J{G}(u)$ has
nonnegative integer coefficients, but no combinatorial interpretation
of these integers is known; we shall return to this point as well
in Section 6.

\section{The Hermite-Biehler Theorem}

The Hermite-Biehler Theorem is a very useful criterion which determines
whether a polynomial with real coefficients has all its zeros in the left
half-plane.  For a nonzero $P(u)\in\RR[u]$, if every $u\in\CC$ such that
$P(u)=0$ satisfies  $\Re(u)\leq 0$ then $P(u)$ is {\em Hurwitz
quasi-stable}.  (For the relevance of this concept to solutions of
linear ordinary differential equations, see Theorem 3.1 of Barnett \cite{Ba}.)

Suppose that $A,B\in\RR[x]$ both have only real zeros, that
those of $A$ are $\xi_1 \leq ... \leq \xi_a$ and that those of $B$ are
$\theta_1 \leq ... \leq \theta_b$.  We say that $A$ {\em
interlaces} $B$ if $\deg B=1+\deg A$ and the zeros of $A$ and $B$ satisfy
$\theta_1\leq\xi_1\leq\theta_2\leq\cdots\leq\xi_a\leq\theta_{a+1}.$
We also say that $A$ {\em alternates left of} $B$ if $\deg A=\deg B$ and
the zeros of $A$ and $B$ satisfy
$\xi_1\leq\theta_1\leq\xi_2\leq\cdots\leq\xi_a\leq\theta_a.$
We use the notation $A\prec B$ for ``either $A$ interlaces $B$ or
$A$ alternates left of $B$.''  (Any polynomial which
stands in this  relation {\em a fortiori} has only real zeros.)
This is a {\em closed condition} in the sense that if $A_{n}$ and
$B_{n}$ are sequences of polynomials converging to $A$ and $B$,
respectively, and if $A_{n}\prec B_{n}$ for all $n\geq 0$, then $A\prec B$.
By convention we say that for any polynomial $A$ with only
real zeros, both $A\prec 0$ and $0\prec A$ hold.  A polynomial
is {\em standard} when either it is identically zero or its
leading coefficient is positive. For brevity, we say that a
polynomial {\em has only nonpositive zeros} to indicate that
either it is identically zero or all of its zeros are real and
nonpositive.  Henceforth,  if we omit the argument of a polynomial
then we intend that it is a function of the variable $x=u^{2}$.

\begin{THM}[Hermite-Biehler]  Let $P(u)=P_{0}(u^2)+uP_{1}(u^2)
\in\RR[u]$ be standard.  Then $P(u)$ is Hurwitz quasi-stable if
and only if both $P_{0}$ and $P_{1}$ are standard, have only
nonpositive zeros, and $P_{1}\prec P_{0}$. \end{THM}

The proof of the Hermite-Biehler Theorem in Gantmacher \cite{Ga} covers
only the case of polynomials for which all zeros have strictly negative
real part, but the statement given here can be deduced from it easily
by a limiting argument.  The following lemmas will be useful; Lemmas 2.2
and 2.3 can be proven using the techniques from Section 3 of \cite{W1}
and Lemma 2.4 can be proven using the techniques from Section 5 of \cite{W1}.

\begin{LMA} Let $P_{1},\ldots,P_{k}$ be polynomials in
$\RR[x]$, all of which have only real zeros and none of which
is identically zero.  If $P_{1}\prec P_{2}\prec\cdots\prec P_{k}$ and
$P_{1}\prec P_{k}$ then $P_{i}\prec P_{j}$ for all $1\leq i\leq j\leq k$.
\end{LMA}

\begin{LMA}  Let $A,P,Q$ be standard polynomials in $\RR[x]$ which
have only nonpositive zeros, and assume that $A\not\equiv 0$.\\
{\rm (a)}\ \ $P\prec Q$ if and only if $Q\prec xP$.\\
{\rm (b)}\ \ If $A\prec P$ and $A\prec Q$ then $A\prec P+Q$.\\
{\rm (c)}\ \ If $P\prec A$ and $Q\prec A$ then $P+Q\prec A$.\\
{\rm (d)}\ \ If $P\prec Q$ then $P\prec P+Q\prec Q$.
\end{LMA}

\begin{LMA}  Let $P,Q$ be standard polynomials in $\RR[x]$ which
have only nonpositive zeros.  Then $P\prec Q$ if and only if for
all $\lam,\rho> 0$, both $\lam P+\rho Q$ and
$\lam Q+\rho xP$ have only nonpositive zeros.
\end{LMA}

\section{Thick Cacti}

Returning to the case of a network $G$, we define polynomials
$\Ji{G}$ and $\Jo{G}$ in $\NN[x]$ by separating $\J{G}(u)$
into its odd and even parts, respectively:
\showon \J{G}(u)=\Jo{G}(u^{2})+u\Ji{G}(u^{2}).\showoff
From (3.1) and the Hermite-Biehler Theorem it follows that (1.8) and
(1.14) are each equivalent to  \showon \Ji{G}\prec\Jo{G}.\showoff
From (1.15) and (3.1) it follows that if $G$ and $N$ are networks with
exactly one vertex in common then
\showon \left\{\begin{array}{rcrcr}
\Jo{G\cup N} &=& \Jo{G}\Jo{N} &+& x\Ji{G}\Ji{N},\\
\Ji{G\cup N} &=& \Jo{G}\Ji{N} &+& \Ji{G}\Jo{N}.
\end{array}\right.\showoff
For each natural number $c$ we define $E_{c}$ and $O_{c}$ in $\NN[x]$
by \showon (u+1)^{c}=E_{c}(u^{2})+uO_{c}(u^{2}).\showoff
From (1.16), (3.1), and (3.4) it follows that if $\sigma$ is a
$c$-spindle in the network $G$ then for $c$ even
\showon \left\{\begin{array}{rcrcrcr}
\Jo{G} &=& E_{c}\Jo{G\drop\sigma} &+& xO_{c}\Ji{G\drop\sigma}
&+& xO_{c}\Ji{G/\sigma},\\
\Ji{G} &=& E_{c}\Ji{G\drop\sigma} &+& O_{c}\Jo{G\drop\sigma}
&+& O_{c}\Jo{G/\sigma},
\end{array}\right.\showoff
and for $c$ odd \showon \left\{\begin{array}{rcrcrcr}
\Jo{G} &=& E_{c}\Jo{G\drop\sigma} &+& xO_{c}\Ji{G\drop\sigma}
&+& E_{c}\Jo{G/\sigma},\\
\Ji{G} &=& E_{c}\Ji{G\drop\sigma} &+& O_{c}\Jo{G\drop\sigma}
&+& E_{c}\Ji{G/\sigma}.
\end{array}\right.\showoff
With $G$, $G^{\bul}$, $N$, and $N^{\bul}$ as in (1.6), (1.11), and (1.17),
let
$\J{}(u):=\J{G}(u)$ and $\J{\tri}(u):=\J{G}(u)+\J{G^{\bul}}(u)$
and $K(u):=\J{N}(u)$ and $K_{\tri}(u):=\J{N}(u)+\J{N^{\bul}}(u)$;
we find that
\showon \left\{\begin{array}{rcrcrcrcr}
\Jo{G\cup N} &=& \Jo{}\Ko{\tri} &+& \Jo{\tri}\Ko{} &+&
                 x\Ji{}\Ki{\tri} &+& x\Ji{\tri}\Ki{},\\
\Ji{G\cup N} &=& \Jo{}\Ki{\tri} &+& \Jo{\tri}\Ki{} &+&
                 \Ji{}\Ko{\tri} &+& \Ji{\tri}\Ko{}. \end{array}\right.\showoff

\begin{LMA}  For any natural numbers $a$ and $b$,
$O_{a+b}\prec E_{a+b}$,
$E_{a}O_{b} \prec E_{a+b}$, $O_{a}E_{b}\prec E_{a+b}$,
and $O_{a}O_{b} \prec O_{a+b} \prec E_{a}E_{b}$.
\end{LMA}
\begin{proof}  We proceed by induction on $a+b$, the base $a+b\leq 1$ being
evident.  If $a=0$ then the only nontrivial claim is that
$O_{b}\prec E_{b}$, since $O_{0}=0$ and $E_{0}=1$; we can
prove this claim by considering the case $a':=1$ and $b':=b-1$
instead.  Similarly, we may also dispense with the case $b=0$,
so assume that $a\geq 1$ and $b\geq 1$.
From $(u+1)^{a+b}=(u+1)^{a}(u+1)^{b}$ it follows that
\showon\left\{\begin{array}{lcrcr}
E_{a+b} &=& E_{a}E_{b} &+& xO_{a}O_{b},\\
O_{a+b} &=& E_{a}O_{b} &+& O_{a}E_{b}.\end{array}\right.\showoff
By induction, we have $O_{a}\prec E_{a}$ and $O_{b}\prec E_{b}$,
and hence $O_{a}O_{b}\prec E_{a}O_{b}\prec E_{a}E_{b}$ and
$O_{a}O_{b}\prec O_{a}E_{b}\prec E_{a}E_{b}$.  Lemma 2.3 now
implies the result.  \end{proof}

A {\em cactus} is a connected simple graph in which each edge is
contained in at most one cycle.  We now consider the special case of
networks $G$ such that $G^{\natural}$ is a cactus.  It is convenient
to introduce the notations, for each natural number $c$,
\showon \eps(c):=\left\{\begin{array}{ll}
1 & \ \ \mbox{if $c$ is even},\\
0 & \ \ \mbox{if $c$ is odd},\end{array}\right.
\ \ \ \mbox{and}\ \ \ \delta(c):=1-\eps(c),\showoff
and \showon S_{c}:=\left\{\begin{array}{ll}
O_{c} & \ \ \mbox{if $c$ is even},\\
E_{c} & \ \ \mbox{if $c$ is odd}.\end{array}\right.\showoff
Denoting a network with two vertices and $m$ edges by $mT_{2}$ we
have for all positive integers $m$,
\showon J_{\eps(m)}^{mT_{2}}=S_{m} \ \ \mbox{and}\ \
J_{\delta(m)}^{mT_{2}}\equiv 0.\showoff
For a finite list of positive integers $\mathbf{c}:=(c_{1},\ldots,c_{n-1})$,
let $T[\mathbf{c}]$ denote any network $T$ such that $T^{\natural}$ is a
tree and the sizes of the spindles of $T$ are given by the list
$\mathbf{c}$\ (so $n_{T}=n$).  From (3.3), (3.11), and induction on $n$ it
follows that
\showon J_{\eps(n+m)}^{T[\mathbf{c}]}=x^{\nu[\mathbf{c}]}S_{c_{1}}\cdots S_{c_{n-1}}
\ \ \mbox{and}\ \ J_{\delta(n+m)}^{T[\mathbf{c}]}\equiv 0,\showoff
where $m=c_{1}+\cdots+c_{n-1}$ and
$\nu[\mathbf{c}]:=\lfloor(\eps(c_{1})+\cdots+\eps(c_{n-1}))/2\rfloor$.
The condition $\Ji{G}\prec\Jo{G}$ is trivial only in this case,
as the following proposition shows.

\begin{PROP}  If $G$ is a network such that either $\Ji{G}\equiv 0$
or $\Jo{G}\equiv 0$ then $G^{\natural}$ is a tree.
\end{PROP}
\begin{proof}  If $G^{\natural}$ is not a tree then let $\sigma$ be a spindle
of $G$ corresponding to an edge of $G^{\natural}$ which is not a
cut-edge.  Since $G\drop\sigma$ is connected, $\R{G\drop\sigma}(q)
\not\equiv 0$, so $\J{G\drop\sigma}(u)\not\equiv 0$, so either
$\Ji{G\drop\sigma}\not\equiv 0$ or $\Jo{G\drop\sigma}\not\equiv 0$.
Since all $J$-polynomials of networks have nonnegative
coefficients, (3.5) and (3.6) imply that $\Ji{G}\not\equiv 0$ and
$\Jo{G}\not\equiv 0$. \end{proof}

Lemmas 2.3 and 3.1, formulas (3.5), (3.6), and (3.12), and induction
on $n$ suffice to prove Theorem 3.3;  we omit the details since
we obtain the most interesting part of Theorem 3.3 (the condition
$\Ji{C[\mathbf{c}]}\prec\Jo{C[\mathbf{c}]}$) as a special case of Corollary 5.3.
The crucial simplifying feature in the proof of Theorem 3.3 is that,
in the notation of (3.12), $J_{\delta(n+m)}^{T[\mathbf{c}]}\equiv 0$.

\begin{THM}  Let $\mathbf{c}:=(c_{1},\ldots,c_{n})$ be a  sequence of
$n\geq 3$ positive integers, and let $C[\mathbf{c}]$ denote any network
with spindles of sizes $c_{1},\ldots,c_{n}$ such that the
underlying simple graph is a cycle.  Then $J_{\delta(n+m)}^{C[\mathbf{c}]}
=(n-1)x^{\nu[\mathbf{c}]}S_{c_{1}}\cdots S_{c_{n}}$ and  $\Ji{C[\mathbf{c}]}
\prec\Jo{C[\mathbf{c}]}$, where $m=c_{1}+\cdots+c_{n}$.\end{THM}

\begin{CORO}  Let $G$ be a network for which $G^\natural$ is a cactus.
Then $R_{G}(q)$ is Schur quasi-stable.\end{CORO}
\begin{proof}  By (1.8) and (1.9) it suffices to prove the result for two-connected
networks.  A two-connected network $G$ satisfying the hypothesis
is such that $G^\natural$ is either an edge or a cycle.  The result
follows from (3.2), (3.11), and Theorem 3.3. \end{proof}

\section{Interpolatory Hypercubes of Polynomials}

An {\em interpolatory $0$-cube} is a standard polynomial $A$ which has
only nonpositive zeros.  An {\em interpolatory $1$-cube} is a pair
$(A,B)$ of standard polynomials which have only nonpositive zeros
and are such that $A\prec B$.  We present the theory next for interpolatory
squares and then generalize to higher-dimensional hypercubes.  The
starting point is an analogue of the ``Box Lemma,'' Theorem 5.4 of \cite{W1}.

\begin{PROP} Let $A,B,P,Q$ be in $\RR[x]$,
and consider the following two conditions:\\
$C_{1}:$\ \ For any $\lam,\rho> 0:\ \lam A+\rho B
\prec \lam P+\rho Q$ and $\lam B+\rho xA \prec
\lam Q+\rho xP$ are interpolatory $1$-cubes.\\
$C_{2}:$\ \ For any $\kappa,\pi> 0:\ \kappa A+\pi P
\prec \kappa B+\pi Q$ and $\kappa P+\pi xA \prec
\kappa Q+\pi xB$ are interpolatory $1$-cubes.\\
These  conditions $C_{1}$ and $C_{2}$ are equivalent. \end{PROP}
\begin{proof} It follows from Lemma 2.4 that each of the conditions
$C_{1}$ and $C_{2}$ is equivalent to the condition that for all
$\kappa,\lam,\pi,\rho>0$, each of the polynomials
$\kappa\lam A+\kappa\rho B+\pi\lam P+\pi\rho Q$,
$\kappa\lam B+\kappa\rho xA+\pi\lam Q+\pi\rho xP$,
$\kappa\lam P+\kappa\rho Q+\pi\lam  xA+\pi\rho xB$, and
$\kappa\lam Q+\kappa\rho xP+\pi\lam xB+\pi\rho x^{2}A$
have only nonpositive zeros. \end{proof}

\noindent
Notice that if the conditions of Proposition 4.1 hold then, since
$\prec$ is a closed condition,  in fact they hold for all
$\kappa,\lam,\pi,\rho\geq 0$.

In the diagrams which follow it is convenient to use an arrow
$A\rightarrow B$ to denote $A\prec B$.    If the equivalent conditions
of Proposition 4.1 hold for the  polynomials $A,B,P,Q$ we say that
$$ \begin{array}{ccc}
B & \lraw & Q\\
\uaw & \SQ & \uaw\\
A & \lraw & P\end{array}$$
is an {\em interpolatory square}, and use the notation $\boxplus$
to indicate this.  Notice that if one of the three squares
\showon \begin{array}{ccc}B & \lraw & Q\\
\uaw & & \uaw\\ A & \lraw & P\end{array}
\ \ \mbox{or}\ \
\begin{array}{ccc} xA & \lraw & xP\\
\uaw & & \uaw\\ B & \lraw & Q\end{array}
\ \ \mbox{or}\ \
\begin{array}{ccc} Q& \lraw & xB\\
\uaw & & \uaw\\ P & \lraw & xA\end{array}\showoff
is interpolatory then all three are, by Lemma 2.3.

\begin{LMA} Let $A,B,P,Q$ be in $\RR[x]$.
If $A\prec B$ and $P\prec Q$ are interpolatory $1$-cubes then
$$ \begin{array}{ccc} BP & \lraw & BQ\\
\uaw & \SQ & \uaw\\ AP & \lraw & AQ\end{array}.$$
\end{LMA}
\begin{proof}  Condition $C_{1}$  of Proposition 4.1 is verified easily
by using Lemma 2.3. \end{proof}

\begin{LMA}  Consider $A,B,P,Q,S,T$ in $\RR[x]$,
with either $A\not\equiv 0$ or $B\not\equiv 0$.
$$\mbox{If}\ \ \begin{array}{ccccc}
T & \llaw & B & \lraw & Q\\
\uaw & \SQ  & \uaw & \SQ & \uaw\\
S & \llaw & A & \lraw & P\end{array}
\ \ \mbox{then}\ \
\begin{array}{ccc}
B & \lraw & Q+T\\
\uaw & \SQ & \uaw\\
A & \lraw & P+S\end{array}.$$
\end{LMA}
\begin{proof} For any $\lam,\rho>0$ we have
$\lam A+\rho B\prec\lam P+\rho Q$ and
$\lam A+\rho B\prec\lam S+\rho T$ since the squares are
interpolatory.  Thus,
$\lam A+\rho B\prec \lam(P+S)+\rho(Q+T)$ by Lemma 2.3.
Also,
$\lam B+\rho xA\prec\lam Q+\rho xP$ and
$\lam B+\rho xA\prec\lam T+\rho xS$ since the squares are
interpolatory.  Thus,
$\lam B+\rho xA\prec \lam(Q+T)+\rho x(P+S)$ by Lemma 2.3.
We have verified condition $C_{1}$ of Proposition 4.1, and hence
the result.\end{proof}

\noindent
It follows from Lemma 4.3 and (4.1) that, under the hypothesis of
Lemma 4.3,
\showon\mbox{if}\ \ \begin{array}{ccccc}
T & \lraw & B & \llaw & Q\\
\uaw & \SQ  & \uaw & \SQ & \uaw\\
S & \lraw & A & \llaw & P\end{array}
\ \ \mbox{then}\ \
\begin{array}{ccc}
Q+T & \lraw & B\\
\uaw & \SQ & \uaw\\
P+S & \lraw & A\end{array}.\showoff

\begin{LMA}  Let $A,B$ be in $\RR[x]$.
If $A\prec B$ is an interpolatory $1$-cube then
$$\begin{array}{ccc} B & \lraw & xA\\
\uaw & \SQ & \uaw\\ A & \lraw & B\end{array}.$$\end{LMA}
\begin{proof}  If either $A\equiv 0$ or $B\equiv 0$ then the result is
trivial, so assume that $A\not\equiv 0$ and $B\not\equiv 0$.
For any $\lam,\rho>0$, Lemma 2.3 implies that
$A\prec\lam A+\rho B\prec B\prec \lam B+\rho xA\prec xA,$
and since $A\prec xA$ it follows from Lemma 2.2 that $\lam A+
\rho B\prec \lam B+\rho xA$.  Also by Lemma 2.3, the condition
that $\lam B +\rho xA \prec \lam xA+\rho xB$ is equivalent to
the condition that $\lam A+\rho B\prec \lam B+\rho xA$, which
we have just shown, and so condition $C_{1}$ of Proposition 4.1 is
verified. \end{proof}

Using Lemmas 4.2, 4.3, and 4.4, one may adapt the proof of Lemma
3.1 to show that for all $a\geq 0$ and $b\geq 0$:
\showon\begin{array}{ccccc}
E_{a}O_{b} & \lraw & E_{a+b} & \lraw & xO_{a}E_{b}\\
\uaw & \SQ  & \uaw & \SQ & \uaw\\
O_{a}O_{b} & \lraw & O_{a+b} & \lraw & E_{a}E_{b}\end{array}.\showoff

To extend these ideas from squares to hypercubes of
any dimension we must first introduce some notation.
Fix an integer $k\geq 0$, and let $\PP:(\ZZ/2\ZZ)^{k}
\rightarrow\RR[x]$.
If $k\geq 1$ then for $\lam,\rho\geq 0$ we let
$_{1}\I_{\lam}^{\rho}$ denote the operator which maps $\PP$ to
the $2^{k-1}$ polynomials
$\Q:(\ZZ/2\ZZ)^{k-1}\rightarrow\RR[x]$ given by $Q_{\al}:=
\lam P_{1\al}+\rho P_{0\al}$ for all $\al\in(\ZZ/2\ZZ)^{k-1}$;
for each $1\leq i\leq k$ we define the {\em $i$-th interpolation
operator}  $_{i}\I_{\lam}^{\rho}$
by a similar interpolation on the $i$-th coordinate of $\PP$.
Given $k$-tuples $\lam:=(\lam_{1},\ldots,\lam_{k})$ and
$\rho:=(\rho_{1},\ldots,\rho_{k})$ of nonnegative real numbers,
we define $\II_{\lam}^{\rho}:=\: _{1}\I_{\lam_{k}}^{\rho_{k}}\cdots\:
_{1}\I_{\lam_{2}}^{\rho_{2}}\: _{1}\I_{\lam_{1}}^{\rho_{1}}$.
For each $1\leq i\leq k$ let $\eta(i):=0\ldots 010\ldots 0$
(with the $1$ in the $i$-th coordinate) be the coordinate vectors
of $(\ZZ/2\ZZ)^{k}$, and denote by $\Phi_{i}$ the {\em $i$-th
flip operator} which associates to $\PP:(\ZZ/2\ZZ)^{k}\rightarrow
\RR[x]$ the $2^{k}$ polynomials
$$(\Phi_{i}\PP)_{\al}:=\left\{\begin{array}{ll}
xP_{\al+\eta(i)} & \ \mbox{if}\ \al_{i}=0,\\
P_{\al+\eta(i)} & \ \mbox{if}\ \al_{i}=1,\end{array}\right.$$
for each $\al\in(\ZZ/2\ZZ)^{k}$.  For any $S\subseteq\{1,\ldots,k\}$
we let $\Phi_{S}:=\prod_{i\in S}\Phi_{i}$.
We say that $\PP:(\ZZ/2\ZZ)^{k}\rightarrow\RR[x]$ is an
{\em interpolatory $k$-cube} of polynomials when the following
condition holds:  for all $k$-tuples $\lam:=(\lam_{1},\ldots,\lam_{k})$
and $\rho:=(\rho_{1},\ldots,\rho_{k})$ of positive real numbers,
and for all $S\subseteq\{1,\ldots,k\}$, the polynomial
$\II_{\lam}^{\rho}\Phi_{S}\PP$ is standard and has only nonpositive zeros.

\begin{PROP}  Fix $k\geq 1$, and let $\PP:(\ZZ/2\ZZ)^{k}\rightarrow\RR[x]$.
Consider the following conditions $C_{i}$ for each $1\leq i\leq k$:\\
$C_{i}:$\ \ For all $\lam,\rho>0$:\ both $_{i}\I_{\lam}^{\rho}\PP$ and
$_{i}\I_{\lam}^{\rho}\Phi_{i}\PP$ are interpolatory $(k-1)$-cubes.\\
The conditions $C_{i}$ for $1\leq i\leq k$ are each equivalent to the
condition that $\PP$ is an interpolatory $k$-cube.
\end{PROP}
\begin{proof} This follows from Lemma 2.4 (the case $k=1$)  as in
the proof of Proposition 4.1 (the case $k=2$). \end{proof}

\begin{LMA}  Fix nonnegative integers $k$ and $\ell$, and let
$\PP:(\ZZ/2\ZZ)^{k}\rightarrow\RR[x]$ and
$\Q:(\ZZ/2\ZZ)^{\ell}\rightarrow\RR[x]$.
Define $\mathbf{S}:(\ZZ/2\ZZ)^{k+\ell}\rightarrow\RR[x]$
by $S_{\al\beta}:=P_{\al}Q_{\beta}$ for all
$\al\in(\ZZ/2\ZZ)^{k}$ and $\beta\in(\ZZ/2\ZZ)^{\ell}$.
If both $\PP$ and $\Q$ are interpolatory hypercubes then $\mathbf{S}$
is an interpolatory $(k+\ell)$-cube.\end{LMA}
\begin{proof} If $k=0$ then it may be checked directly that $\mathbf{S}$
satisfies the definition of an interpolatory $\ell$-cube, and so
we proceed by induction on $k\geq 1$.  For any $\lambda,\rho>0$ we
have, by Proposition 4.5, interpolatory $(k-1)$-cubes
$_{1}\I_{\lam}^{\rho}\PP$ and $_{1}\I_{\lam}^{\rho}\Phi_{1}\PP$.
By induction, both $_{1}\I_{\lam}^{\rho}\mathbf{S}$ and
$_{1}\I_{\lam}^{\rho}\Phi_{1}\mathbf{S}$ are interpolatory
$(k-1+\ell)$-cubes, and so Proposition 4.5 implies that $\mathbf{S}$
is an interpolatory $(k+\ell)$-cube. \end{proof}

\begin{LMA}  Fix $k\geq 1$, and let $\PP,\Q:(\ZZ/2\ZZ)^{k}\rightarrow\RR[x]$
be such that $P_{1\al}=Q_{1\al}$ for all
$\al\in(\ZZ/2\ZZ)^{k-1}$, and $P_{1\al}\not\equiv 0$ for at least one
$\al\in(\ZZ/2\ZZ)^{k-1}$.  Define $\mathbf{S}:(\ZZ/2\ZZ)^{k}\rightarrow\RR[x]$ by
$$S_{\al}:=\left\{\begin{array}{ll}
P_{\al}+Q_{\al}& \ \ \mbox{if}\ \ \al_{1}=0,\\
P_{\al} & \ \ \mbox{if}\ \ \al_{1}=1,\end{array}\right.$$
for all $\al\in(\ZZ/2\ZZ)^{k}$.  If both $\PP$ and $\Q$ are interpolatory
$k$-cubes then $\mathbf{S}$ is an interpolatory $k$-cube.\end{LMA}
\begin{proof}  Again we argue by induction on $k$, the base $k=1$ being
Lemma 2.3(b), and the case $k=2$ being Lemma 4.3;  for the induction
step let $k\geq 2$.  For any $\lam,\rho>0$:\
$_{2}\I_{\lam}^{\rho}\PP$ and $_{2}\I_{\lam}^{\rho}\Q$
are interpolatory $(k-1)$-cubes satisfying the hypothesis,
and hence $_{2}\I_{\lam}^{\rho}\mathbf{S}$
is an interpolatory $(k-1)$-cube.  Also,  for any $\lam,\rho>0$:\
$_{2}\I_{\lam}^{\rho}\Phi_{2}\PP$ and
$_{2}\I_{\lam}^{\rho}\Phi_{2}\Q$ are interpolatory
$(k-1)$-cubes satisfying the hypothesis, and hence
$_{2}\I_{\lam}^{\rho}\Phi_{2}\mathbf{S}$ is an interpolatory $(k-1)$-cube.
This verifies condition $C_{2}$ of Proposition 4.5 for $\mathbf{S}$.\end{proof}

\noindent
By an argument analogous to the derivation of (4.2) from Lemma 4.3, we
may conjugate everything in Lemma 4.7 by $\Phi_{1}$ to obtain another
similar statement (which is left to the reader).

\begin{LMA} Fix $k\geq 0$, and let $\PP:(\ZZ/2\ZZ)^{k}
\rightarrow\RR[x]$.  Define $\Q:(\ZZ/2\ZZ)^{k+1}
\rightarrow\RR[x]$ by $Q_{1\al}:=P_{\al}$ and $Q_{0\al}:=
(\Phi_{1}\PP)_{\al}$ for all $\al\in(\ZZ/2\ZZ)^{k}$.  If $\PP$
is an interpolatory $k$-cube then $\Q$ is an interpolatory
$(k+1)$-cube.\end{LMA}
\begin{proof} The case $k=0$ is obvious, and the case $k=1$ is Lemma 4.4;
for $k\geq 2$ we proceed by induction on $k$.  Choose $\lam,\rho>0$ and
apply the induction hypothesis to the interpolatory $(k-1)$-cubes
$_{2}\I_{\lam}^{\rho}\PP$ and $_{2}\I_{\lam}^{\rho}\Phi_{2}\PP$ to find that
$_{3}\I_{\lam}^{\rho}\Q$ and $_{3}\I_{\lam}^{\rho}\Phi_{3}\Q$ are
interpolatory $k$-cubes.  By Proposition 4.5 it follows that $\Q$ is an
interpolatory $(k+1)$-cube. \end{proof}

\begin{LMA}  Fix $k\geq 2$, and let
$\PP:(\ZZ/2\ZZ)^{k}\rightarrow\RR[x]$.  Define
$\Q:(\ZZ/2\ZZ)^{k-1}\rightarrow\RR[x]$ by
$$Q_{1\al}:=P_{01\al}+P_{10\al}\ \ \mbox{and}\ \
Q_{0\al}:=P_{00\al}+xP_{11\al}$$
for all $\al\in(\ZZ/2\ZZ)^{k-2}$.  Assume that there is a
$\beta\in(\ZZ/2\ZZ)^{k-2}$ such that either $P_{00\beta}\not\equiv 0$ or
$P_{11\beta}\not\equiv 0$, and that there is a $\gamma\in(\ZZ/2\ZZ)^{k-2}$
such that either $P_{01\gamma}\not\equiv 0$ or $P_{10\gamma}\not\equiv 0$.
If $\PP$ is an interpolatory $k$-cube then $\Q$ is an interpolatory
$(k-1)$-cube.\end{LMA}
\begin{proof} Since $\PP$ is an interpolatory $k$-cube, each of
$_{1}\I_{0}^{1}\PP$, $_{1}\I_{1}^{0}\Phi_{2}\PP$,
$_{2}\I_{0}^{1}\PP$, and $_{2}\I_{1}^{0}\Phi_{1}\PP$ is an
interpolatory $(k-1)$-cube.  By exchanging the first and second coordinates
if necessary, we may assume that $\PP$ is such that $P_{10\gamma}\not\equiv 0$
for some $\gamma\in(\ZZ/2\ZZ)^{k-2}$.  Thus we may apply Lemma 4.7
to $_{1}\I_{1}^{0}\Phi_{2}\PP$ and $_{2}\I_{0}^{1}\PP$;  denote the result by
$\mathbf{S}$.  If $P_{01\gamma}\equiv 0$ for all $\gamma\in(\ZZ/2\ZZ)^{k-2}$ then
$\mathbf{S}=\Q$, which suffices to prove the result.  Otherwise, we may also
apply Lemma 4.7 to $_{1}\I_{0}^{1}\PP$ and $_{2}\I_{1}^{0}\Phi_{1}\PP$;  denote
the result by $\mathbf{T}$.  Then $_1\I_0^1\mathbf{S}=\: _1\I_0^1\mathbf{T}=\:
_1\I_0^1 \Q$, and by the hypothesis there is some $\beta\in(\ZZ/2\ZZ)^{k-2}$
such that $Q_{0\beta}\not\equiv 0$.  Thus we may apply the $\Phi_1$-conjugate
form of Lemma 4.7 to $\mathbf{S}$ and $\mathbf{T}$;  the result is $\Q$, which
proves the result. \end{proof}

Of course, by conjugating with a permutation of indices one may apply Lemma
4.9 to any two coordinates $1\leq i<j\leq k$ of $\PP$.  In this case the
correspondence between the indices of $\PP$ and the indices of
$\Q$ will be taken to be $\ell\mapsto\ell$ for $1\leq\ell<j$,
$j\mapsto i$, and $\ell\mapsto\ell-1$ for $j<\ell\leq k$, generalizing
the case of $i=1$ and $j=2$ in the statement above.

\section{Series-Parallel Networks}

After much experimentation one arrives at the following hypothesis.
For a network $G$ and distinct vertices $v$ and $w$ of $G$, let
$G^-$ be obtained from $G$ by deleting all edges between $v$ and
$w$, and let $G^{\bul}$ be obtained from $G$ by identifying
$v$ and $w$ and removing any loops thus created.  We shall say that
$\{v,w\}$ is {\em very amicable in $G$} if

\showon \begin{array}{ccc}
\Jo{\bul} & \lraw & x\Ji{-}\\
\uaw & \SQ & \uaw\\
\Ji{\bul} & \lraw & \Jo{-}\end{array},
\ \ \mbox{and equivalently}\ \
\begin{array}{ccc}
\Jo{-} & \lraw & x\Ji{\bul}\\
\uaw & \SQ & \uaw\\
\Ji{-} & \lraw & \Jo{\bul}\end{array},
\showoff
where $J_-(u):=\J{G^-}(u)$ and $\J{\bul}(u):=\J{G^{\bul}}(u)$.
(These conditions are equivalent, by (4.1).)   In fact, this condition is
too strong, and we shall say that $\{v,w\}$ is {\em amicable
in $G$} if the condition
\showon \begin{array}{ccc}
\Jo{-}+\Jo{\bul} & \lraw & x\Ji{-}\\
\uaw & \SQ & \uaw\\
\Ji{-}+\Ji{\bul} & \lraw & \Jo{-}\end{array}\showoff
is satisfied.  Notice that $J_{-}(u)\equiv 0$ if and only if $v$ and
$w$ are adjacent in $G$ and $v\sim w$ is a cut-edge of
$G^{\natural}$.  In this case, (5.1) and (5.2) are each equivalent to
$\Ji{\bul}\prec\Jo{\bul}$; otherwise, from (5.1) we have $\Ji{-}\prec\Jo{-}$,
to which we apply Lemma 4.4, and then (4.2) and (5.1) imply (5.2).
In either case, if $\{v,w\}$ is very amicable in $G$ then $\{v,w\}$ is
amicable in $G$.  Notice that $J_{\bul}(u)\not\equiv 0$ since $G$, and
hence $G^{\bul}$, is connected.

\begin{LMA} Let $G$ be a network and let $\{v,w\}$ be amicable in $G$.
Then $\Ji{G}\prec\Jo{G}$.\end{LMA}
\begin{proof} Let $G^-$, $G^{\bul}$, $J_-(u)$, and $J_{\bul}(u)$ be as in the
above paragraph, and let $J_{\tri}(u):=J_-(u)+J_{\bul}(u)$.
Let $v$ and $w$ be joined by exactly $a$ edges of $G$. The hypothesis
that $\{v,w\}$ is amicable in $G$ is (5.2).
Applying Lemma 4.6 to (5.2) and $O_a\prec E_a$ gives an
interpolatory $3$-cube.
\showon \begin{array}{ccccccc}
O_a\Jo{\tri} & \lraw & xO_a\Ji{-}& &
E_a \Jo{\tri} & \lraw & xE_a\Ji{-}\\
\uaw & \SQ & \uaw & \Longrightarrow & \uaw & \SQ & \uaw\\
O_a\Ji{\tri} & \lraw & O_a\Jo{-}& &
E_a\Ji{\tri} & \lraw & E_a\Jo{-}\end{array}\showoff
Index the coordinates of (5.3) by $1,2,3$ in the order $\uaw,\rightarrow,
\Rightarrow$.  If $a>0$ and $J_-(u)\not\equiv 0$ then we may apply Lemma
4.9 to coordinates $2$ and $3$ of (5.3), yielding an interpolatory square.
\showon \begin{array}{ccc}
E_a\Jo{\tri}+xO_a\Ji{-} & \lraw & xE_a\Ji{-}+xO_a\Jo{\tri}\\
\uaw & \SQ & \uaw\\
E_a\Ji{\tri}+O_a\Jo{-} & \lraw & E_a\Jo{-}+xO_a\Ji{\tri}
\end{array}\showoff
If $a=0$ then $O_a=0$, and (5.4) is obtained from (5.3) by applying
$_3\I_0^1$;  if $J_-(u)\equiv 0$ then (5.4) is obtained from (5.3) by
applying $_2\I_1^0\Phi_3$.  Thus, in all cases (5.4) is an interpolatory
square.  If $a$ is odd then (3.6) and the left column of (5.4) give
$\Ji{G}\prec\Jo{G}$, while if $a$ is even then (3.5) and the right column
of (5.4) give $\Jo{G}\prec x\Ji{G}$.  \end{proof}

\begin{THM}
Let $G$ and $N$ be networks which intersect in exactly one
vertex $v$, let $w'\neq v$ be a vertex of $G$, and let $w''\neq v$
be a vertex of $N$.  Let $U:=G\cup N$, let $W$ denote the network obtained
from $U$ by identifying $w'$ and $w''$, and let $w$ denote the
image of $w'$ and $w''$ in $W$.  If $\{v,w'\}$ is amicable in $G$
and $\{v,w''\}$ is amicable in $N$ then:\\
{\rm (a)}\ \ $\{v,w'\}$ and $\{v,w''\}$ are amicable in $U$, and\\
{\rm (b)}\ \ $\{v,w\}$ is amicable in $W$, and\\
{\rm (c)}\ \ $\{w',w''\}$ is very amicable in $U$.
\end{THM}
\begin{proof}
Let $v$ and $w'$ be joined by exactly $a$ edges of $G$, and let
$v$ and $w''$ be joined by exactly $b$ edges of $N$.  Let $G^-$ denote
the network obtained by deleting the $a$ edges between $v$ and $w'$ in $G$,
and let $G^{\bul}$ denote the network obtained by identifying $v$ and $w'$
in $G$ and removing the $a$ loops thus produced.   Let $N^-$ denote
the network obtained by deleting the $b$ edges between $v$ and $w''$ in $N$,
and let $N^{\bul}$ denote the network obtained by identifying $v$ and $w''$
in $N$ and removing the $b$ loops thus produced.  Let $W^-$ denote
the network obtained by deleting the $a+b$ edges between $v$ and $w$ in $W$,
and let $W^{\bul}$ denote the network obtained by identifying $v$ and $w$
in $W$ and removing the $a+b$ loops thus produced.
To simplify notation, let $J(u):=\J{G}(u)$, $J_-(u):=\J{G^-}(u)$, and
$\J{\bul}(u):=\J{G^{\bul}}(u)$, let $K(u):=\J{N}(u)$, $K_-(u):=\J{N^-}(u)$,
and $K_{\bul}(u):= \J{N^{\bul}}(u)$, and let $L(u):=\J{W}(u)$,
$L_-(u):=\J{W^-}(u)$, and $L_{\bul}(u):= \J{W^{\bul}}(u)$.  We will also
use the notations $\J{\tri}(u):=\J{-}(u)+\J{\bul}(u)$,
$K_{\tri}(u):=K_-(u)+K_{\bul}(u)$, and $L_{\tri}(u):=L_-(u)+L_{\bul}(u)$.

By Lemma 5.1 we see that $\Ji{}\prec\Jo{}$ and $\Ki{}\prec\Ko{}$.
Now apply Lemma 4.6 to  (5.2) (with $w'$ in place of $w$)
and $\Ki{}\prec\Ko{}$ to get an interpolatory $3$-cube.
\showon \begin{array}{ccccccc}
\Jo{\tri}\Ki{} & \lraw & x\Ji{-}\Ki{}& &
\Jo{\tri}\Ko{} & \lraw & x\Ji{-}\Ko{}\\
\uaw & \SQ & \uaw & \Longrightarrow & \uaw & \SQ & \uaw\\
\Ji{\tri}\Ki{} & \lraw & \Jo{-}\Ki{}& &
\Ji{\tri}\Ko{} & \lraw & \Jo{-}\Ko{}\end{array}\showoff
Index the coordinates of (5.5) as for (5.3).  If $\Ki{}\not\equiv 0$ and
$\Ko{}\not\equiv 0$ then, since $J_{\bul}\not\equiv 0$, we may apply Lemma
4.9 to coordinates $1$ and $3$ of (5.5), yielding the interpolatory square
\showon \begin{array}{ccc}
\Jo{\tri}\Ko{}+x\Ji{\tri}\Ki{} & \lraw & x\Jo{-}\Ki{}+x\Ji{-}\Ko{}\\
\uaw & \SQ & \uaw\\
\Jo{\tri}\Ki{}+\Ji{\tri}\Ko{}& \lraw & \Jo{-}\Ko{}+x\Ji{-}\Ki{}
\end{array},\showoff
which from (3.3) is seen to be
\showon\begin{array}{ccc}
\Jo{G^-\cup N}+\Jo{G^{\bul}\cup N} & \lraw & x\Ji{G^-\cup N}\\
\uaw & \SQ & \uaw\\
\Ji{G^-\cup N}+\Ji{G^{\bul}\cup N} & \lraw & \Jo{G^-\cup N}\end{array}.\showoff
If $\Ki{}\equiv 0$ then (5.6) is obtained from (5.5) by applying
$_3\I_{0}^{1}$;  if $\Ko{}\equiv 0$ then (5.6) is obtained from (5.5) by
applying $_{3}\I_{1}^{0}\Phi_{2}$.  In all cases (5.7) is an interpolatory
square, showing that $\{v,w'\}$ is amicable in $U$.  Since the hypothesis
is symmetric in $G$ and $N$ we also conclude that $\{v,w''\}$ is
amicable in $U$, proving part (a).

For part (b), apply Lemma 4.6 to (5.2) and its analogue for $N$ to get an
interpolatory $4$-cube.
\showon \begin{array}{ccccccc}
\Jo{\tri}\Ko{\tri} & \lraw & x\Ji{-}\Ko{\tri} &  &
x\Jo{\tri}\Ki{-} & \lraw & x^2\Ji{-}\Ki{-}\\
\uaw & \SQ  & \uaw & \Longrightarrow & \uaw & \SQ & \uaw\\
\Ji{\tri}\Ko{\tri} & \lraw & \Jo{-}\Ko{\tri} &  &
x\Ji{\tri}\Ki{-} & \lraw & x\Jo{-}\Ki{-}\\
 & \Uparrow & & \SQ & & \Uparrow & \\
\Jo{\tri}\Ki{\tri} & \lraw & x\Ji{-}\Ki{\tri} &  &
\Jo{\tri}\Ko{-} & \lraw & x\Ji{-}\Ko{-}\\
\uaw & \SQ & \uaw & \Longrightarrow & \uaw & \SQ  & \uaw\\
\Ji{\tri}\Ki{\tri} & \lraw & \Jo{-}\Ki{\tri} &  &
\Ji{\tri}\Ko{-} & \lraw & \Jo{-}\Ko{-}
\end{array}\showoff
Index the coordinates of (5.8) by $1,2,3,4$ in the order $\uaw$,
$\rightarrow$, $\Uparrow$, $\Rightarrow$.  If $J_-(u)\not\equiv 0$
and $K_-(u)\not\equiv 0$ then we may apply Lemma 4.9 to coordinates
$2$ and $4$ of (5.8) to obtain an interpolatory $3$-cube.
\showon \begin{array}{ccc}
x\Jo{\tri}\Ki{-}+x\Ji{-}\Ko{\tri} & \lraw & x\Jo{\tri}\Ko{\tri}+x^2\Ji{-}\Ki{-}\\
\uaw & \SQ & \uaw\\
\Jo{-}\Ko{\tri}+x\Ji{\tri}\Ki{-} & \lraw & x\Jo{-}\Ki{-}+x\Ji{\tri}\Ko{\tri}\\
& \Uparrow &\\
\Jo{\tri}\Ko{-}+x\Ji{-}\Ki{\tri} & \lraw & x\Jo{\tri}\Ki{\tri}+x\Ji{-}\Ko{-}\\
\uaw & \SQ & \uaw\\
\Jo{-}\Ki{\tri}+\Ji{\tri}\Ko{-} & \lraw & \Jo{-}\Ko{-}+x\Ji{\tri}\Ki{\tri}
\end{array}\showoff
If $J_-(u)\equiv 0$ then (5.9) is obtained from (5.8) be applying
$_2\I_1^0\Phi_4$ (and permuting coordinates); if $K_-(u)\equiv 0$ then
(5.9) is obtained from (5.8) by applying $_4\I_1^0\Phi_2$.  In all cases,
(5.9) is an interpolatory $3$-cube.
If $J_-(u)\not\equiv 0$ or $K_-(u)\not\equiv 0$ then we may apply Lemma
4.9 to coordinates $1$ and $3$ of (5.9), yielding
\showon \begin{array}{ccc}
\begin{array}{r} x\Jo{-}\Ki{\tri}+x\Ji{-}\Ko{\tri}\\
+x\Jo{\tri}\Ki{-}+x\Ji{\tri}\Ko{-} \end{array} & \lraw &
\begin{array}{r} x\Jo{-}\Ko{-}+x^{2}\Ji{-}\Ki{-}\\
+x\Jo{\tri}\Ko{\tri}+x^2\Ji{\tri}\Ki{\tri} \end{array} \\
\uaw & \SQ & \uaw\\
\begin{array}{r} \Jo{-}\Ko{\tri}+x\Ji{-}\Ki{\tri}\\
+\Jo{\tri}\Ko{-}+x\Ji{\tri}\Ki{-}  \end{array} & \lraw &
\begin{array}{r} x\Jo{-}\Ki{-}+ x\Ji{-}\Ko{-}\\
+x\Jo{\tri}\Ki{\tri}+x\Ji{\tri}\Ko{\tri} \end{array}
\end{array}\showoff
If $J_-(u)\equiv 0$ and $K_-(u)\equiv 0$ then first assume that
neither $(G^\bul)^\natural$ nor $(N^\bul)^\natural$ is a tree.
By Proposition 3.2 we may apply Lemma 4.9 to (5.9) to produce (5.10).
Otherwise, if $\Ji{\bul}\equiv 0$ then (5.10) is obtained from (5.9)
by applying $_2\I_0^1$, if $\Jo{\bul}\equiv 0$ then (5.10) is obtained
from (5.9) by applying $_2\I_1^0\Phi_3$, and similarly in case $\Ko{\bul}
\equiv 0$ or $\Ki{\bul}\equiv 0$.  In all cases, (5.10) is an interpolatory
square.  From (1.17) we have
$L_-(u)=J_-(u)K_{\tri}(u)+J_{\tri}(u)K_-(u)$, and hence
$L_{\tri}(u)=L_-(u)+J_{\bul}(u)K_{\bul}(u)=J_-(u)K_-(u)+J_{\tri}(u)K_{\tri}(u)$.
From this one sees that (5.10) is
\showon \begin{array}{ccc}
xL_1^- & \lraw & xL_0^{\tri} \\
\uaw & \SQ & \uaw \\
L_0^- & \lraw & xL_1^{\tri}\end{array},\showoff
and from (4.1) it follows that $\{v,w\}$ is amicable in $W$,
proving part (b).

For part (c) we begin with (5.4) and its analogue for $N$, that is
\showon \begin{array}{ccc}
E_b\Ko{\tri}+xO_b\Ki{-} & \lraw & xE_b\Ki{-}+xO_b\Ko{\tri}\\
\uaw & \SQ & \uaw\\
E_b\Ki{\tri}+O_b\Ko{-} & \lraw & E_b\Ko{-}+xO_b\Ki{\tri}
\end{array}.\showoff
As in the proof of Lemma 5.1, both (5.4) and (5.12) are interpolatory 
squares, so that by Lemma 4.6 we obtain an interpolatory $4$-cube
$\Q$; we index the coordinates of $\Q$ so that $1$ and $3$
correspond to $\uaw$ and $\rightarrow$ in (5.4) and $2$ and $4$
correspond to $\uaw$ and $\rightarrow$ in (5.12), respectively.
The cases in which either $a>0$ and $G^{\natural}$ is a tree or
$b>0$ and $N^{\natural}$ is a tree are slightly degenerate;  assume
first that neither condition holds.  Then we can apply Lemma 4.9 to
coordinates $1$ and $2$ of $\Q$ to obtain an interpolatory $3$-cube
$\mathbf{T}$;  the entries of $\mathbf{T}$ are as follows.

\begin{eqnarray*}
T_{000}
&=& (xE_{a}\Ji{-}+xO_{a}\Jo{\tri})(xE_{b}\Ki{-}+xO_{b}\Ko{\tri})\\
& & +x(E_{a}\Jo{-}+xO_{a}\Ji{\tri})(E_{b}\Ko{-}+xO_{b}\Ki{\tri})
\end{eqnarray*}
\begin{eqnarray*}
T_{001}
&=& (xE_{a}\Ji{-}+xO_{a}\Jo{\tri})(E_{b}\Ko{\tri}+xO_{b}\Ki{-})\\
& & +x(E_{a}\Jo{-}+xO_{a}\Ji{\tri})(E_{b}\Ki{\tri}+O_{b}\Ko{-})
\end{eqnarray*}
\begin{eqnarray*}
T_{010}
&=& (E_{a}\Jo{\tri}+xO_{a}\Ji{-})(xE_{b}\Ki{-}+xO_{b}\Ko{\tri})\\
& & +x(E_{a}\Ji{\tri}+O_{a}\Jo{-})(E_{b}\Ko{-}+xO_{b}\Ki{\tri})
\end{eqnarray*}
\showon\begin{array}{lcl}
T_{011}
&=& (E_{a}\Jo{\tri}+xO_{a}\Ji{-})(E_{b}\Ko{\tri}+xO_{b}\Ki{-})\\
& & +x(E_{a}\Ji{\tri}+O_{a}\Jo{-})(E_{b}\Ki{\tri}+O_{b}\Ko{-})\\ \\
T_{100}
&=& (xE_{a}\Ji{-}+xO_{a}\Jo{\tri})(E_{b}\Ko{-}+xO_{b}\Ki{\tri})\\
& & +(E_{a}\Jo{-}+xO_{a}\Ji{\tri})(xE_{b}\Ki{-}+xO_{b}\Ko{\tri})
\end{array}\showoff
\begin{eqnarray*}
T_{101}
&=& (xE_{a}\Ji{-}+xO_{a}\Jo{\tri})(E_{b}\Ki{\tri}+O_{b}\Ko{-})\\
& & +(E_{a}\Jo{-}+xO_{a}\Ji{\tri})(E_{b}\Ko{\tri}+xO_{b}\Ki{-})
\end{eqnarray*}
\begin{eqnarray*}
T_{110}
&=& (E_{a}\Jo{\tri}+xO_{a}\Ji{-})(E_{b}\Ko{-}+xO_{b}\Ki{\tri})\\
& & +(E_{a}\Ji{\tri}+O_{a}\Jo{-})(xE_{b}\Ki{-}+xO_{b}\Ko{\tri})
\end{eqnarray*}
\begin{eqnarray*}
T_{111}
&=& (E_{a}\Jo{\tri}+xO_{a}\Ji{-})(E_{b}\Ki{\tri}+O_{b}\Ko{-})\\
& & +(E_{a}\Ji{\tri}+O_{a}\Jo{-})(E_{b}\Ko{\tri}+xO_{b}\Ki{-})
\end{eqnarray*}
If $a>0$ and $G^\natural$ is a tree then $J_-(u)\equiv 0$ and either
$\Ji{\bul}\equiv 0$ or $\Jo{\bul}\equiv 0$;  if $\Ji{\bul}\equiv 0$
then $\mathbf{T}=\: _{1}\I_{0}^{1}\Q$, while if $\Jo{\bul}\equiv 0$ then
$\mathbf{T}=\: _{1}\I_{1}^{0}\Phi_{1}\Q$.  The case when $b>0$ and $N^\natural$
is a tree is handled similarly.  In all cases $\mathbf{T}$ is an
interpolatory $3$-cube.  Notice that $J_{\bul}(u)\not\equiv 0$ and
if $a=0$ then $J_{-}(u)\not\equiv 0$, and similarly for $N$;  from this
it follows that $T_{0\al}+T_{1\al}\not\equiv 0$ for all $\al\in(\ZZ/2\ZZ)^{2}$.

If $a$ and $b$ are both odd then $\Ji{}=E_a\Ji{\tri}+O_a\Jo{-}$ and
$\Jo{}=E_a\Jo{\tri}+xO_a\Ji{-}$ and $\Ki{}=E_b\Ki{\tri}+O_b\Ko{-}$ and
$\Ko{}=E_b\Ko{\tri}+xO_b\Ki{-}$ and $L_1=E_{a+b}L_1^-+O_{a+b}L_0^{\tri}$
and $L_0=E_{a+b}L_0^-+xO_{a+b}L_1^{\tri}$.  Thus we find that
$T_{011}=\Jo{}\Ko{}+x\Ji{}\Ki{}$ and $T_{111}=\Jo{}\Ki{}+\Ji{}\Ko{}$
and, by using (3.7) and (3.8), that $T_{001}+T_{010}=xL_1$ and $T_{101}+T_{110}=L_0$.
Applying Lemma 4.3 to $_2\I_1^0\mathbf{T}$ and $_3\I_1^0\mathbf{T}$ we see that
\showon \begin{array}{ccc}
T_{011} & \lraw & T_{001}+T_{010}\\
\uaw & \SQ & \uaw\\
T_{111}& \lraw & T_{101}+T_{110}\end{array},\showoff
which shows that $\{w',w''\}$ is very amicable in $U$ in this case.

If $a$ and $b$ are both even then
$\Ji{}=E_a\Ji{-}+O_a\Jo{\tri}$ and $\Jo{}=E_a\Jo{-}+xO_a\Ji{\tri}$ and
$\Ki{}=E_b\Ki{-}+O_b\Ko{\tri}$ and $\Ko{}=E_b\Ko{-}+xO_b\Ki{\tri}$ and
$L_1=E_{a+b}L_1^{-}+O_{a+b}L_0^{\tri}$ and $L_0=E_{a+b}L_0^{-}+
xO_{a+b}L_1^{\tri}$.  Thus we find that $T_{000}=x(\Jo{}\Ko{}+x\Ji{}\Ki{})$
and $T_{100}=x(\Jo{}\Ki{}+\Ji{}\Ko{})$ and, by using (3.7) and (3.8), that
$T_{001}+T_{010}=xL_1$ and
$T_{101}+T_{110}=L_0$.  Applying (4.2) to $_2\I_0^1\mathbf{T}$ and
$_3\I_0^1\mathbf{T}$ we see that
\showon \begin{array}{ccc}
T_{001}+T_{010} & \lraw & T_{000}\\
\uaw & \SQ & \uaw\\
T_{101}+T_{110} & \lraw & T_{100}\end{array},\showoff
which shows that $\{w',w''\}$ is very amicable in $U$ in this case,
by (4.1).

In the remaining case, $a$ and $b$ have opposite parity; by symmetry we may
assume that $a$ is even and $b$ is odd.  Thus
$\Ji{}=E_a\Ji{-}+O_a\Jo{\tri}$ and $\Jo{}=E_a\Jo{-}+xO_a\Ji{\tri}$ and
$\Ki{}=E_b\Ki{\tri}+O_b\Ko{-}$ and $\Ko{}=E_b\Ko{\tri}+xO_b\Ki{-}$ and
$L_1=E_{a+b}L_1^{\tri}+O_{a+b}L_0^-$ and $L_0=E_{a+b}L_0^{\tri}+xO_{a+b}L_1^-$.
Thus we find that $T_{001}=x(\Jo{}\Ki{}+\Ji{}\Ko{})$ and
$T_{101}=\Jo{}\Ko{}+x\Ji{}\Ki{}$ and, by using (3.7) and (3.8), that
$T_{000}+xT_{011}=xL_0$ and
$T_{100}+xT_{111}=xL_1$.  Applying Lemma 4.3 to $_2\I_1^0\mathbf{T}$ and
$_3\I_1^0\Phi_2\mathbf{T}$ we see that
\showon \begin{array}{ccc}
T_{001} & \lraw & T_{000}+xT_{011}\\
\uaw & \SQ & \uaw\\
T_{101}& \lraw & T_{100}+xT_{111}\end{array},\showoff
which shows that $\{w',w''\}$ is very amicable in $U$ in this case,
by (4.1).  This completes the proof. \end{proof}

\begin{CORO}  If $G$ is a series-parallel network with terminals
$\{s,t\}$ then $\{s,t\}$ is amicable in $G$, and hence
$\Ji{G}\prec\Jo{G}$.\end{CORO}
\begin{proof} The basis of induction $m=1$ is clear, by (3.11).
For the induction step, let $G'$ be a series parallel-network with $m\geq 2$;
so $G'$ can be written either as a series connection or as a parallel
connection of series-parallel networks $G$ and $N$ both with strictly
fewer edges than $G'$.  If the connection is series then Theorem 5.2(c)
provides the induction step; if the connection is parallel then Theorem 
5.2(b)
provides the induction step. Lemma 5.1 then completes the proof. \end{proof}

Theorem 0.2 now follows immediately, since (1.9) shows that a network
satisfies (1.8) if and only if each of its two-connected components does,
and Corollary 5.3 proves that (3.2), and hence (1.8), holds for all networks
in the class $\mathfrak{SP}$.

In fact, the argument proving Theorem 5.2(a) can be used to show that if
$G$ and $N$ intersect in exactly one vertex $v$, if $\{w_{1},w_{2}\}$ is
amicable in $G$, and if $\Ji{N}\prec\Jo{N}$, then $\{w_{1},w_{2}\}$ is
amicable in $G\cup N$.  One can use this and Theorem 5.2 to show that if
$G^{\natural}$ is a cactus then every pair $\{v,w\}$ of distinct vertices
of $G$ is amicable in $G$.  By Lemma 5.1, the following conjecture
implies Conjecture 0.1.

\begin{CONJ}  Let $G$ be a network and let $v\neq w$ be
vertices of $G$.  Then $\{v,w\}$ is amicable in $G$.\end{CONJ}

\noindent In fact, I believe that Conjecture 5.4 is just the
first in a hierarchy of conditions on $G$ involving interpolatory
hypercubes of arbitrary dimension.  Determining what these
conditions are might lead to an inductive proof of them all, and
hence of Conjecture 0.1, using the technique of Section 4.

\section{$f$-Vectors of Matroids}

Let $E$ be a set with $m$ elements, and let $\Omega$ be a
{\em set system} (or ``hypergraph'') on $E$, that is, $\Omega$
is any collection of subsets of $E$; members of $\Omega$ will be
called {\em faces} of $\Omega$. Let $d_{\Omega}:=
\max\{\#S:\ S\in\Omega\}$ be the {\em degree} of $\Omega$, and
for $0\leq i\leq d$ let $f_{i}(\Omega)$ be the number of faces
$S\in\Omega$ such that $\#S=i$;  we define the {\em rank-generating
function} of $\Omega$ to be
\showon F_{\Omega}(z)=f_{0}+f_{1}z+f_{2}z^{2}+\cdots+f_{d}z^{d}.\showoff
We may factor this polynomial as $F_{\Omega}(z)=(1+z)^{d-t}
\widetilde{F}_{\Omega}(z)$ such that $\widetilde{F}_{\Omega}(-1)
\neq 0$, defining the {\em subdegree} $t_{\Omega}:=\deg
\widetilde{F}_{\Omega}(z)$ of $\Omega$, and the coefficients of
$\widetilde{F}_{\Omega}(z)=\sum_{i=0}^{t}\tilde{f}_{i}z^{i}$
in the process.  Clearly
\showon \tilde{f}_{i}=\sum_{\ell=0}^{\infty}\binom{d-t+\ell-1}{\ell}
(-1)^{\ell}f_{i-\ell}\ \ \mbox{for all}\ \ 0\leq i\leq t,\showoff
and
\showon f_{i}=\sum_{\ell=0}^{d-t}\binom{d-t}{\ell}\tilde{f}_{i-\ell}
\ \ \mbox{for all}\ \ 0\leq i\leq d,\showoff
with the conventions that $f_{i}=\tilde{f}_{i}=0$ if $i<0$ and
$f_{i}=0$ if $i>d$ and $\tilde{f}_{i}=0$ if $i>t$.

The {\em reliability function} $R_{\Omega}(q)$ of $\Omega$ is the
probability that, if each element of $E$ is selected independently
with probability $0\leq q\leq 1$, then the random subset $\mathcal{S}(q)
\subseteq E$ consisting of the selected elements of $E$ is a face of
$\Omega$. By partitioning the event that $\mathcal{S}(q)\in\Omega$ into
its constituent subevents one sees immediately
that \showon R_{\Omega}(q)=\sum_{i=0}^{d}f_{i}q^{i}(1-q)^{m-i}.\showoff
Thus, we can write $R_{\Omega}(q)=(1-q)^{m-d}H_{\Omega}(q)$,
where the {\em $H$-polynomial} of $\Omega$ is defined by
\showon H_{\Omega}(q):=(1-q)^{d}F_{\Omega}\left(\frac{q}{1-q}\right).\showoff
A simple calculation shows that $\deg H_{\Omega}(q)=t_{\Omega}$ and
that $H_{\Omega}(q)=\sum_{i=0}^{t}h_{i}q^{i}$ depends only upon
$\widetilde{F}_{\Omega}(z)$.

Certain classes of set systems are of special interest with respect
to these polynomials. Let $\mathfrak{S}$ denote the class of
simplicial complexes, let $\mathfrak{M}$ denote the class of
(simplicial complexes of independent sets of) matroids, let
$\mathfrak{G}^{*}$ denote the class of cographic matroids, and let
$\mathfrak{BC}$ denote the class of set systems $\Omega$ for
which $H_{\Omega}(q)$ is Schur quasi-stable.  For a network $G$,
let $M:=M^{*}(G)$ be the cographic matroid associated with $G$; then
the polynomials $R_{G}(q)$, $H_{G}(q)$, and $J_{G}(u)$ defined in
Sections $0$ and $1$ equal the polynomials $R_{M}(q)$, $H_{M}(q)$, and
$J_{M}(u)$ defined in this section, respectively.   Thus, the
Brown-Colbourn Conjecture is that $\mathfrak{G}^{*}$ is a subclass of
$\mathfrak{BC}$.

Lemma 6.1 and Proposition 6.2 were suggested by Theorem 4.3 and the
remark on page 585 of Brown and Colbourn \cite{BC}.  For a positive
integer $k$ and a set system $\Omega$ on the set $E$, we define
the set system $k\Omega$ on the set $E\times\{1,\ldots,k\}$ as
follows:  $\{(e_{1},i_{1}),\ldots,(e_{r},i_{r})\}\subseteq
E\times\{1,\ldots,k\}$ is a face of $k\Omega$ if and only if
$\{e_{1},\ldots,e_{r}\}$ are pairwise distinct elements of $E$,
and this set is a face of $\Omega$.

\begin{LMA} Let $\Omega$ be a set system on a set
$E$ of size $m$ and let $k$ be a positive integer.  Then
$$\R{k\Omega}(q)=((1-q)^{k}+kq(1-q)^{k-1})^{m}
\R{\Omega}\left(\frac{kq}{1+(k-1)q}\right).$$\end{LMA}
\begin{proof} For each of the $m$ elements $e\in E$, at most one of the
elements $(e,1),\ldots,(e,k)$ can be selected if the random subset
of selected elements $\mathcal{S}(q)$ is to be a face of $k\Omega$;
these events occur independently, each with probability
$(1-q)^{k}+kq(1-q)^{k-1}$.  Conditioning on the conjunction of these
events, the conditional probability that exactly one of $(e,1),
\ldots,(e,k)$ is selected is
\showon\hat{q}:=\frac{kq(1-q)^{k-1}}{(1-q)^{k}+kq(1-q)^{k-1}}
=\frac{kq}{1+(k-1)q}\showoff
for each $e\in E$, and hence the conditional probability that
$\mathcal{S}(q)\in\Omega$ is $R_{\Omega}(\hat{q})$. \end{proof}

\begin{PROP} For any set system $\Omega$, there is an integer
$K_{\Omega}$ such that for all $k\geq K_{\Omega}$, $k\Omega$ is in
the class $\mathfrak{BC}$.
\end{PROP}
\begin{proof} Let $\Omega$ be defined on a set $E$ with $m$ elements.
With $\hat{q}$ defined as in (6.6) we have
$R_{k\Omega}(q)=(1-q)^{km-m}(1+(k-1)q)^{m}R_{\Omega}(\hat{q}).$
The zeros of $R_{k\Omega}(q)$ due to the factors
$(1-q)^{km-m}(1+(k-1)q)^{m}$ are inside the unit disc $|q|\leq 1$
for all $k\geq 1$.  If $\xi\in\CC$ is such that $R_{\Omega}(\xi)=0$
then each factor $(\hat{q}-\xi)$ of $R_{\Omega}(\hat{q})$
contributes a zero of $R_{k\Omega}(q)$ at $q_{0}:=\xi/(\xi+k(1-\xi))$.
If $\xi=1$ then $q_{0}=1$, and if $\xi\neq 1$ then we can choose
$k$ sufficiently large that $|q_{0}|<1$.  Since $R_{\Omega}(q)$ has
only finitely many zeros, there is some $K_{\Omega}$ such that
$k\geq K_{\Omega}$ suffices for all factors, proving the result.
\end{proof}

\noindent
In fact, the proof of Brown and Colbourn \cite{BC} shows that if
$M$ is a matroid then $K_{M}=d_{M}+1$ suffices in Proposition 6.2,
although they do not state this explicitly.

Proposition 6.3 provides some weak support for the idea
that all matroids are in the class $\mathfrak{BC}$, but at present
there is not enough evidence to really justify any opinion on this
strengthening of the Brown-Colbourn Conjecture.

\begin{PROP}  For $1\leq d<m$, let $U_{m}^{d}$ denote the uniform
matroid of rank $d$ with $m$ elements, let $F_{m}^{d}(z)$ be
its rank-generating function, and construct $H_{m}^{d}(q)$ as
in $(6.5)$.  If $q\in\CC$ is such that $H_{m}^{d}(q)=0$ then
$(m-d)^{-1}\leq|q|\leq d(m-1)^{-1}$.  In particular, $U_{m}^{d}$ is in the class
$\mathfrak{BC}$.\end{PROP}
\begin{proof}  For all $1\leq d<m$ we have $F_{m}^{d}(z)=\sum_{i=0}^{d}
\binom{m}{i} z^{i}$, and so $F_{m}^{d}(-1)\neq 0$.
For $d=1$ this gives $F_{m}^{1}(z)=1+mz$ and
$H_{m}^{1}(q)=1+(m-1)q$, satisfying the statement of the proposition.
From the familiar recurrence relations for binomial coefficients
it follows that for all $1<d<m$, $F_{m}^{d}(z)=zF_{m-1}^{d-1}(z)+
F_{m-1}^{d}(z)$ and $H_{m}^{d}(q)=qH_{m-1}^{d-1}(q)+H_{m-1}^{d}(q)$.
By induction, one sees that for all $1\leq d<m$,
$$H_{m}^{d}(q)=\sum_{i=0}^{d}\binom{m-d-1+i}{i} q^{i}.$$
The successive ratios of these coefficients are $\lam_i:=\binom{m-d-1+i}{i}
\binom{m-d+i}{i+1}^{-1}=(i+1)(m-d+i)^{-1}$, which are
nondecreasing as $i$ runs from $0$ to $d-1$.  Thus, by the
Enestr\"om-Kakeya Theorem (see Theorem B of Anderson, Saff, and Varga
\cite{ASV}) it follows that all complex zeros of
$H_{m}^{d}(q)$ satisfy $\lam_0\leq|q|\leq \lam_{d-1}$. \end{proof}

The {\em $\tilde{f}$-vector} $(\tilde{f}_{0},\ldots,\tilde{f}_{t})$
of a set system $\Omega$ in the class $\mathfrak{BC}$ must satisfy some
strong inequalities, as we now explain; when $t_{\Omega}=d_{\Omega}$
this $\tilde{f}$-vector agrees with the {\em $f$-vector}
$(f_{0},\ldots,f_{d})$ of $\Omega$.  We introduce the {\em $J$-polynomial}
of a set system $\Omega$ by defining
\showon
\quad\quad\J{\Omega}(u):=(-2)^{t}\widetilde{F}_{\Omega}
\left(\frac{-1-u}{2}\right),
\ \ \ \mbox{so that}\ \ \
\widetilde{F}_{\Omega}(z)=\frac{\J{\Omega}(-1-2z)}{(-2)^t},\showoff
where we have the relations $u:=-1-2z$ and $z=(-1-u)/2$.
In terms of the coefficients $\J{\Omega}(u)=
\sum_{k=0}^{t}j_{k}u^{k}$ this relation is
\showon j_k=\sum_{i=k}^t\binom{i}{k}(-2)^{t-i}\tilde{f}_i
\ \ \mbox{for all}\ \ 0\leq k\leq t,\showoff
and conversely,
\showon \tilde{f}_i=2^{i-t}\sum_{k=i}^t\binom{k}{i}(-1)^{t-k}j_{k}
\ \ \mbox{for all}\ \ 0\leq i\leq t.\showoff
The relation between $J_{\Omega}(u)$ and $H_{\Omega}(q)$ is as in
(1.12) and (1.13).
By reasoning analogous to that showing the equivalence of (1.8) and (1.14),
one sees that a set system $\Omega$ is in the class
$\mathfrak{BC}$ if and only if
$\J{\Omega}(u)$ is Hurwitz quasi-stable.  A theorem of Asner
\cite{A} (see also Kemperman \cite{K}) states that a polynomial
$J(u)=\sum_{k=0}^{t}j_{k}u^{k}$ in $\RR[u]$ with $j_{t}>0$
is such that all of its zeros have strictly negative real part
if and only if every minor of the {\em Hurwitz matrix}
\showon
\mathbf{H}[J(u)]:=\left[\begin{array}{cccccc}
j_{0} & 0 & 0 & \cdots & 0 & 0\\
j_{2} & j_{1} & j_{0} & \cdots & \ddots & 0\\
j_{4} & j_{3} & j_{2} & \ddots & \vdots & \vdots\\
\vdots & \vdots & \ddots & j_{t-2} & j_{t-3} & j_{t-4}\\
0 & \ddots & \cdots & j_{t} & j_{t-1} & j_{t-2}\\
0 & 0 & \cdots & 0 & 0 & j_{t}\end{array}\right]\showoff
is nonnegative, and $\det\mathbf{H}[J(u)]>0$.  One direction of this
equivalence survives in the limit (the other does not): if $J(u)$
is Hurwitz quasi-stable then every minor of $\mathbf{H}[J(u)]$ is
nonnegative (see \cite{A,K}).
We let $\mathfrak{BC}'$ denote the class of set systems $\Omega$
such that every minor of $\mathbf{H}[\J{\Omega}(u)]$ is nonnegative;
this class contains $\mathfrak{BC}$ (and hence, by Theorem 0.2, the
cographic matroid of each network in $\mathfrak{SP}'$).  Also, we denote
by $\mathfrak{J}_{+}$ the class of set systems $\Omega$ such that
$j_{k}(\Omega)\geq 0$ for all $0\leq k\leq t_{\Omega}$; this class
contains $\mathfrak{BC}'$.

For example, consider the simplicial complex $I$ consisting of the
faces of the icosahedron.  We have $F_{I}(z)=1+12z+30z^{2}+20z^{3}$
and so $t_{I}=d_{I}=3$, and the calculation of $\J{I}(u)$ can be
illustrated by
\showon\begin{array}{r|r|rrrr}
1 & -8  & -8 & & & \\
12 & 4  & 48 & 48 & & \\
30 & -2 & -60 & -120 & -60 & \\
20 & 1  & 20 & 60 & 60 & 20 \\ \hline
& & 0 & -12 & 0 & 20 \end{array}\showoff
so that $\J{I}(u)=-12u+20u^{3}$.  This is an example of a simplicial
polytope which is not in the class $\mathfrak{J}_{+}$.
As another example, let $\Theta$ be the broken-circuit complex
(see Brylawski \cite{Bry}) of (the graphic matroid of) $K_{2,3}$.
Then $F_{\Theta}(z)=(1+z)(1+5z+10z^{2}+7z^{3})$, so $d_{\Theta}=4$
and $t_{\Theta}=3$ and we calculate that $\J{\Theta}(u)=
-1+u+u^{2}+7u^{3}$; this $\Theta$ is a
broken-circuit complex which is not in $\mathfrak{J}_{+}$.
Simplicial polytopes, broken-circuit complexes, and matroids
are each subclasses of the class of Cohen-Macaulay complexes;
see Stanley \cite{St2}.  For a discussion of the location of zeros
of $F_{\Delta}(z)$ for Cohen-Macaulay complexes $\Delta$ in general, see
\cite{W2}.

Our last theorem also supports the possibility that $\mathfrak{M}$ might
be a subclass of $\mathfrak{BC}$.

\begin{THM} Every matroid is in the class $\mathfrak{J}_+$.\end{THM}
\begin{proof}  Let $M$ be a matroid of rank $d$ which has exactly $c$ coloops,
and let $M'$ be the
matroid obtained by deleting all loops and coloops of $M$.
Then $F_{M}(z)=(1+z)^{c}F_{M'}(z)$, so that $\widetilde{F}_{M}(z)
=\widetilde{F}_{M'}(z)$, and $M'$ has no loops or coloops.  Since
$H_{M}(q)$ and $J_{M}(u)$ depend only upon $\widetilde{F}_{M}(z)$,
we may replace $M$ by $M'$ and henceforth assume that $M$ has no
loops or coloops.

It is a standard result of matroid theory (see (7.12) of Bj\"orner
\cite{Bj}, for example) that the Tutte polynomial $T_{M}(x,y)$
of a matroid $M$ may be specialized to yield
\showon T_{M}(x,1)=h_{0}x^{d}+h_{1}x^{d-1}+\cdots+h_{d-1}x+h_{d}=
x^{d}H_{M}(1/x),\showoff
where $h_{i}=0$ if $t_{M}<i\leq d_{M}$. Another standard result
is that if $M$ has no coloops then $h_d>0$; this
follows from Theorem 6.2.13(v) in Brylawski and Oxley
\cite{BO} (and that fact that $b_{ij}\geq 0$ for all $i,j$ in their
notation, see p.127 of \cite{BO}).  (In other words, $t_{M}=d_{M}$
for a matroid with no coloops.)
Chari \cite{Ch1} proves that, since $M$ has no coloops,
there exist integers $s_{i}\geq 1$ for $1\leq i\leq h_{d}$ and
$r_{i\ell}\geq 0$ for $1\leq i\leq h_{d}$ and $1\leq\ell\leq s_{i}$ such that
\showon T_{M}(x,y)=\sum_{i=1}^{h_{d}}\prod_{\ell=1}^{s_{i}}
(y+x+x^{2}+\cdots+x^{r_{i\ell}}).\showoff
Letting $a_{i}:=d-\sum_{\ell}r_{i\ell}$ for $1\leq i\leq h_{d}$ we
obtain from (6.12) and (6.13) that
\showon H_{M}(q)=\sum_{i=1}^{h_{d}}q^{a_{i}}
\prod_{\ell=1}^{s_{i}}(1+q+q^{2}+\cdots+q^{r_{i\ell}}),\showoff
in which each term has degree $d$.  Therefore, applying the
relation (1.12) to (6.14) we obtain
\showon J_{M}(u)=\sum_{i=1}^{h_{d}}(u+1)^{a_{i}}
\prod_{\ell=1}^{s_{i}}\left[\frac{(u+1)^{r_{i\ell}+1}-(u-1)^{r_{i\ell}+1}}{2}
\right],\showoff
which evidently has nonnegative integer coefficients.
Therefore $M$ is in the class $\mathfrak{J}_{+}$. \end{proof}

Theorem 6.4 raises the problem of interpreting the coefficients
of the $J$-polynomial of a matroid combinatorially;  although
one can use (6.15) as a guide, a solution to this problem is not
presently at hand.

The proof of Theorem 0.3 is now clear.  If $M$ is a matroid with
no coloops then $t_{M}=d_{M}$ as in the proof of Theorem 6.4,
and thus the $\tilde{f}$-vector of $M$ coincides with the
$f$-vector of $M$.  By Theorem 6.4, $M$ is in the class
$\mathfrak{J}_{+}$, and hence the conclusion of Theorem 0.3 follows
from (6.8).


\newpage

\begin{thebibliography}{ABC}

\bibitem{ASV} N. Anderson, E.B. Saff, and R.S. Varga,
\textit{On the Enestr\"om-Kakeya theorem and its sharpness}, Linear
Algebra Appl. $\mathbf{28}$ $(1979)$, $5$-$16$.

\bibitem{A} B.A. Asner, Jr., \textit{On the total nonnegativity of the
Hurwitz matrix}, S.I.A.M. J. Appl. Math. $\mathbf{18}$ $(1970)$,
$407$-$414$.

\bibitem{BP} M.O. Ball and J.S. Provan, \textit{Bounds on the
reliability polynomial for shellable independence systems}, S.I.A.M.
J. Alg. Discrete Meth. $\mathbf{3}$ $(1982)$, $166$-$181$.

\bibitem{Ba} S. Barnett, \textit{Polynomials and Linear Control Systems},
Dekker, New York, $1983$.

\bibitem{Bj} A. Bj\"orner, \textit{Homology and shellability of matroids
and geometric lattices}, in:  \textit{Matroid Applications}, (N. White, ed.),
Cambridge University Press, Cambridge, $1992$.

\bibitem{BC} J.I. Brown and C.J. Colbourn, \textit{Roots of the
reliability polynomial}, S.I.A.M. J. Discrete Math. $\mathbf{5}$
$(1992)$, $571$-$585$.

\bibitem{BC2} J.I. Brown and C.J. Colbourn, \textit{Non-Stanley
bounds for network reliability}, J. Algebraic Comb. $\mathbf{5}$
$(1996)$, $13$-$36$.

\bibitem{Bry} T. Brylawski, \textit{The broken-circuit complex},
Trans. Amer. Math. Soc. $\mathbf{234}$ $(1977)$, $417$-$433$.

\bibitem{BO} T. Brylawski and J. Oxley, \textit{The Tutte polynomial
and its applications}, in:  \textit{Matroid Applications}, (N. White, ed.),
Cambridge University Press, Cambridge, $1992$.

\bibitem{Ch1} M. Chari, \textit{Matroid inequalities}, Discrete
Math. $\mathbf{147}$ $(1995)$, $283$-$286$.

\bibitem{Ch2} M. Chari, \textit{Two decompositions in topological
combinatorics with applications to matroid complexes}, Trans.
Amer. Math. Soc., to appear.

\bibitem{Co} C.J. Colbourn, \textit{The Combinatorics of Network
Reliability}, Oxford University Press, New York, $1987$.

\bibitem{Ga} F.R. Gantmacher, \textit{Matrix Theory, vol. II}, Chelsea,
New York, $1960$.

\bibitem{K} J.H.B. Kemperman, \textit{A Hurwitz matrix is totally
positive}, S.I.A.M. J. Math. Anal. $\mathbf{13}$ $(1982)$, $331$-$341$.

\bibitem{St2} R.P. Stanley, \textit{Combinatorics and Commutative Algebra,
Second Edition}, Birkh\"auser, Boston, $1996$.

\bibitem{W1} D.G. Wagner, \textit{Total positivity of Hadamard products},
J. Math. Anal. Appl. $\mathbf{163}$ $(1992)$, $459$-$483$.

\bibitem{W2} D.G. Wagner, \textit{Zeros of rank-generating functions of
Cohen-Macaulay complexes}, Discrete Math. $\mathbf{139}$ $(1995)$, 
$399$-$411$.


\end{thebibliography}
\end{document}